\documentclass[12pt]{amsart}

\usepackage{longtable}
\usepackage{url}
\usepackage{enumerate}
\usepackage{epsfig}
 \usepackage{ifthen}
 \newcommand{\mymarginpar}[1]{%
    \marginpar{\ifthenelse{\isodd{\arabic{page}}}{\flushleft #1}{\flushright #1}}}

 \renewcommand{\phi}{\varphi}
 
 \newcommand{\IC}{\mathbf{C}}                     % komplexe Zahlen
 \newcommand{\IN}{\mathbf{N}}                     % natrliche Zahlen
 \theoremstyle{plain} %%%%%%%%%%%%%%%%%%%%%%%%%%%%%%%%%
 \newtheorem{Theorem}{Theorem}[section]

 \theoremstyle{definition} %%%%%%%%%%%%%%%%%%%%%%%%%%%%
 \newtheorem{Definition}[Theorem]{Definition}

 \numberwithin{equation}{section}

\begin{document}

 \title{Finite von Neumann Algebra Factors with
   Property $\Gamma$}
\dedicatory{Dedicated to Professor R. V. Kadison on his $75^{th}$ birthday} 
 \author{Erik Christensen}
 \address{Institut for Matematiske Fag University of Copenhagen\\
Universitetsparken 5\\
DK 2100 Copenhagen \O \\
Denmark }
 \email{echris@math.ku.dk}
 \date{\today}


 \begin{abstract}
   Techniques introduced by G. Pisier in his proof that finite von
   Neumann factors with property $\Gamma$ have length at most 5 are
   modified to prove that the length is 3. It  is proved that
   if such a factor is a complemented subspace of a larger C*-algebra
   then there exists a projection of norm one from the larger algebra
   onto the smaller and  a  new proof of the fact that the continuous
   Hochschild cohomology group $H^2_c(M,M)$ vanishes  is also included.

 \end{abstract}

 \maketitle

 \section{Introduction}
 \label{sec:introduction}

The articles \cite{CS1,CS2,Pi1,Pi2,Pi3} address among other things 
 the question whether a von Neumann
algebra, on a Hilbert space $H$ which is a complemented subspace of $B(H)$, is
an injective von Neumann algebra, i. e. the image of a projection of
norm 1 from $B(H).$

The question was answered in the positive for all von Neumann algebras
  if the projection onto $M$  is completely
bounded. By an averaging technique one can prove that for properly
infinite von Neumann algebras the existence of a bounded projection
onto $M$ implies the existence of a completely bounded one and hence
show that $ M $ is injective. 

\smallskip

For finite continuous von Neumann algebras this sort of  arguments are not
usable in general so we only have a very scattered or sporadic list
of answers to  this complementation question for finite continuous
factors. The list of algebras for which this  problem has been settled
is identical to the list of algebras to which we can answer the
similarity question  - which we will describe below. It seems that
whenever a method has been developed to answer one of the questions
for a particular algebra, then it more or less immediately yields a
method to construct an answer to the other problem for this particular 
algebra. In this paper we do obtain a certain positive result for the
similarity problem first and then uses this to get a result concerning
algebras which are complemented subspaces.

\smallskip

We will now describe the so called similarity problem \cite{Ka,Pi4}. In
\cite{Ka} Kadison asks whether a homomorphism $\phi$ of a
C*-algebra $A$ into the bounded operators on a Hilbert space $B(H)$ is
similar to a *-homomorphism, meaning that there exists a bounded
invertible operator on $H$ which intertwines the given homomorphism
with a self-adjoint homomorphism. In \cite{Pi4} Pisier lists a lot of
known results and problems related to this question. Over the last
couple of years Pisier has introduced several new ideas 
and obtained remarkable new insights with respect to several
questions related to various forms of similarities. Especially the
concepts called length and similarity degree \cite{Pi5,Pi6,Pi7} 
will play a major role in
this article . We will describe these terms to the
extent to  which we need them in the next section. For now we will just
mention that Pisier has proved that Kadison's question has a positive
answer for a C*-algebra $A$ if and only if $A$ has finite length.
\smallskip

In the article  \cite{Ch2} we proved that Kadison's question could be settled
in the affirmative for finite continuous von Neumann factors with
property $\Gamma$. In modern language we proved that the length
 was at most 44. In \cite{Pi7}  Pisier proves that the length is in the
 interval [3, 5]. Using some basic ideas from \cite{Pi7} and combining them
 with the point of view found in \cite{Ch1} we can actually see that the
 length is 3. This implies that the only computed values ( for C*-algebras )
  actually are 2 and 3,
 and it is quite a mystery if $\infty$ or some other integers will occur as lengths for
 C*-algebras.

 \medskip
 As a bonus we can use the methods to prove that a finite continuous
 factor with property $\Gamma$ which is a complemented subspace of a C*-algebra
  also is
 complemented via a completely positive and completely contractive 
 projection of the C* algebra onto the von Neumann algebra. 

\medskip

If the theory of simultaneously ultra strongly continuous multi linear
mappings were better understood it might be that the methods could
also be used to prove that all the continuous Hochschild cohomology
groups for a factor with property $\Gamma $ vanish. We do have some
results which indicates this, but we have no way of controlling
simultaneous ultra strong continuity except if the cochains are
completely bounded. In the latter case we already know that
that the completely bounded cohomology is in general trivial
\cite{CS2,SS1}. In the case of a 2-cocycle the continuity problems can
be solved quite easily and we have included a new proof of the
vanishing of $H^2_c(M,M)$ for the von Neumann factors with property
$\Gamma$.

\medskip

We refer to the books  \cite{EK,KR} for the theory of operator algebras, 
to \cite{Pa} for results on completely bounded operators, to
\cite{Co,Di,EK} for results on the property $\Gamma$ to
\cite{Pi4,Pi5,Pi6,Pi7} for
results on the similarity question, similarity degree  and on length
and to \cite{SS1} for results on Hochschild cohomology.

 \section{Length and the similarity problem for factors with property $\Gamma$}
 \label{sec:lengthgamma}

The concepts of similarity degree and length was introduced by Pisier
\cite{Pi5,Pi6} as a result of his deep investigations into many
different versions of similarity problems. We will not go into the
details  here but
use the results of these articles for our particular purpose.
 We start by recalling the
definition of length and then describe its deep connections to
the set of bounded homomorphisms of an operator algebra into $B(H)$.

\begin{Definition}
\label{Def:length}
A unital  operator algebra $A$ has finite length at most $l \in \IN $ if there exists
a constant C such that for any $k \in \IN$ and any $x \in
\text{M}_k(A)$ there exist an $n \in \IN$, scalar matrices $\alpha_0
\in \text{M}_{k,n}(\IC),\alpha_1
\in \text{M}_{n}(\IC), \dots ,\alpha_{l-1}
\in \text{M}_{n}(\IC),\alpha_l
\in \text{M}_{n,k}(\IC)$ and diagonal matrices $D_1, \dots , D_l$ in
$\text{M}_n(A)$such that
\begin{displaymath}
x = \alpha_0 D_1 \alpha_1 D_2  \dots D_l
\alpha_l \quad \mathrm{and} \quad (\prod_0^l \parallel \alpha_i \parallel)(\prod_1^l
\parallel D_i \parallel) \leq C\parallel x \parallel.
\end{displaymath}

The length $l(A)$ is defined to be the least  possible  $l$ for which these
conditions are fulfilled. 
\end{Definition}

Suppose $\phi$ is a homomorphism of an operator algebra $A$ into
$B(H)$, then it is quite easy to
see that if $A$ has finite length $l$ with constant $C$ then for any $ k \in \IN $ the
homomorphism $\phi_k : M_k(A) \to M_k(B(H))$ has norm at most $C \|
\phi \|^l $, so  $\phi$ is completely bounded and $ \| \phi \|_{cb}
 \leq C \| \phi \|^l $. In particular finite length implies that any
 bounded  
 homomorphism is completely bounded. The very surprising result
 result of \cite{Pi5}  is that
 the converse is also true, as formulated in the following theorem.

\begin{Theorem}[Pisier]
\label{Thm:simdeg}
Let $ A $ be a unital  operator algebra then any bounded unital 
homomorphism of $A$
into $B(H)$ is completely bounded if and only if
 there exist  positive constants
$C$ and $\alpha$ such  that  for any bounded  unital homomorphism $\phi : A
\to B(K) $; $ \| \phi \|_{cb} \leq C \| \phi \|^{\alpha}$.
Moreover $ A$ has this property if and only if
 $A$ has finite length and the
length $l$ is the minimum over the possible values of $\alpha$
\end{Theorem}

The least $\alpha$ usable above is called the similarity degree and the
theorem tells that the similarity degree and the length is the same integer.
\medskip

In \cite{Pi7} Pisier proves that von Neumann factors with property $\Gamma$
have length at most 5 by constructing a concrete factorization of the type above for certain elements of ``rank'' one or
less in $M_k(A)$. Here we will prove that the length is at most 3
with help of the theorem mentioned above  
by showing that for any bounded homomorphism $\phi $ of a finite
continuous factor with property $\Gamma$ we have $ \| \phi \|_{cb} \leq
 \| \phi \|^3$. We know already by \cite{Ch2,Pi7} that any bounded
 homomorphism of such  a factor
 into $B(H)$ is completely bounded and similar to a
 *-homomorphism. A close examination of the invertible operator which
 intertwines these two homomorphisms - along the lines of the
 computations
in \cite{Pi7} gives the result.

\begin{Theorem}
\label{Thm:gammalength}
Let $ M $ be a continuous finite von Neumann algebra factor with
property $\Gamma$ on a
Hilbert space $H$ then $M$ has length 3.
\end{Theorem}

\begin{proof}

It follows from \cite[Remark 12]{Pi7} that the length is at least 3.   
\smallskip

Let now $\phi$ denote a bounded unital homomorphism of $M$ into $B(K)$
then by \cite{Ch2} there exists  a *-representation $\pi$ of $ M$ on $K$
and an invertible $x$ in $B(K)$ such that $\pi = Ad(x) \phi$. Let 
$h = (x^*x)$ then it follows that $Ad(h^{\frac{1}{2}}) \phi$ also is a *-representation
and we can - and will - replace $x$   by  $h^{\frac{1}{2}}$ and assume  that 
$\pi = Ad(h^{\frac{1}{2}}) \phi$ and $h $ is a positive and invertible
{\em contraction. \/}

\smallskip
 
The weak closure, or the bi-commutant $\pi(M)^{\prime \prime}$  of the unital algebra $\pi(M)$ splits into a sum of an
infinite von Neumann algebra - say $N$ - and a finite one - say
$M_f$.  If one of the summands is missing then the following
computations will all be simplified, so we assume that both $N$ and
$M_f$ are non trivial. The Hilbert space $K$
splits accordingly via a central projection - say $ Q $ - in the weak
closure of $\pi (M)$ such that $N$ acts on $(I-Q)K$ and $M_f$ acts
on $QK$. Finally the representation $\pi$ splits via $Q$ into the sum
of an infinite representation - say $\pi_{\infty}$ and a finite one -
say $\pi_f$. The latter representation is automatically ultra weakly
continuous since any ultra weakly continuous functional on $M_f$ factors
through the trace on $M_f$. On the other hand this trace induces a
trace on $M$ and since this algebra is a factor it has only got a
single ( normalized ) trace. Hence the composition of the
representation $\pi_f$ and an ultra weakly continuous functional on
$M_f$ has a density with respect to the trace on $M$ and therefore
is ultra weakly continuous on $M$. In particular this means that
$\pi_f$ is a normal isomorphism of $M$ onto $M_f$.

\smallskip We will define $c$ by $ c = \| \phi\|^2 $, then for any
unitary $ u$ in $M$ we have $ 0 \leq \phi(u) \phi(u)^* \leq c$.
Since $\phi(u)  =  h^{-\frac{1}{2}}\pi(u)h^{\frac{1}{2}}$
multiplication of 
this inequality from the left and from the right
with $h^{\frac{1}{2}}$ yields, as in \cite[(7)]{Ch1} the following
inequality
 
\begin{equation} \label{eq:1}
\forall u \text{ unitary in } M, \quad \quad  0 \leq \pi(u)h\pi(u)^* \leq ch.
\end{equation}
The unitaries in  M form a group so we can extend the inequality
above to the formally  stronger inequality

\begin{equation} \label{eq:2}
\forall u, v  \text{ unitaries in } M, \quad \quad  0 \leq \pi(u)h\pi(u)^* \leq c\pi(v)h\pi(v)^*.
\end{equation}

Since the Kaplansky Density Theorem makes it possible to approximate unitaries
in the weak closure 
$\pi(M)^{ \prime \prime}$  of $\pi(M)$ strongly with
 unitaries from  $\pi(M)$ the validity of the inequality \ref{eq:2}
 can be extended to

\begin{equation} \label{eq:3}
\forall u, v  \text{ unitaries in } \pi(M)^{ \prime \prime} ,
 \quad \quad  0 \leq uhu^* \leq cvhv^*.
\end{equation}

It is hopefully  clear that the inequality \ref{eq:1} above also imply
that the norm of $\phi$ is at most $c^{\frac{1}{2}}$ and therefore in
order to prove that $M$ has length 3 it is by Theorem  \ref{Thm:simdeg}
sufficient to prove that
for any $n \in \IN$ and any unitary $V = (v_{ij}) \in M_n(M)$ the
inequality just below - named (Goal) - holds.

\begin{equation}
  0 \leq \pi_n(V)(I_{M_n(\IC)} \otimes h) \pi_n(V)^* \leq
    c^3(I_{M_n(\IC)} \otimes  h).\tag{Goal}
\end{equation}

In the rest of the proof we will fix $n \in \IN$ and the unitary 
$V = (v_{ij}) \in M_n(M)$. Further we will use the convention that for any
operator $x \in B(K) $ we will let $\tilde{x} \in M_n(B(K)) $ be given
as $\tilde{x} = I_{M_n(\IC)} \otimes x$.

\smallskip

It is clearly enough to prove the inequality (Goal) for any vector
state on $M_n(B(K))$ so we will also fix a unit vector $\xi = (\xi_1,
\dots , \xi_n) \in \IC^n \otimes K$, a positive real $\varepsilon$ and
verify the inequality (Goal) in the state $\omega_{\xi}$ up to
$\varepsilon$.  

\smallskip We will have to divide the computations
according to the two representations $\pi_{\infty}$ and $\pi_f$. Here
we are faced with the problem that $h$ does not commute with the
central projection $Q$, so we will have to replace $h$ by one which
does. In order to do so we first split $\xi$ as the sum $\xi_{\infty}
+ \xi_f$ by $\xi_{\infty} = \widetilde{(I-Q)}\xi$ and $\xi_f =
\widetilde{Q}\xi$.

Since $M$ is assumed to have property $\Gamma$ we can by Dixmier's
result  \cite[Proposition 1.10]{Di} find a set $\{p_1, \dots , p_n \}$ of pairwise orthogonal
and equivalent projections in $M$ with sum $I$ such that  all the norms
 $\| [v_{ij}, p_l] \|_2$ are small. Here  the norm $\| \,
\|_2$ is the one induced by the pre Hilbert space structure on $M$
coming from the unique trace state. It is a well known fact that
for any ( uniformly ) bounded subset of $M$ the ultra strong topology is
the same as the one coming from this norm. In particular this means that given the vector
$\xi_f$ in $QK \oplus \dots \oplus QK$  and the fact that the
representation $\pi_f$ is normal we can find the set of projections
$\{p_1, \dots , p_n \} \in M$ such that  

\begin{equation} \label{eq:4}
\forall l \in \{1, \dots , n \}: \qquad \|
\widetilde{\pi_f(I-p_l)}(\pi_f)_n(V^*)\widetilde{\pi_f(p_l)}
\xi_f \| \leq 
(\frac{\varepsilon}{n^3})^{\frac{1}{2}}.
 \end{equation}

 Having the projections $\{p_1, \dots , p_n \} \in M$, we will replace
 h by a positive contraction commuting with a finite dimensional
 subfactor of $M_f$ which contains $\{\pi_f(p_1), \dots , \pi_f(p_n)
 \} $ in its main diagonal algebra.  In order to do so we find a set
 of matrix units $(f_{ij})$ in $ M_f$ such that $f_{ii}= \pi_f(p_i)$
 and the set of matrix units generates a subfactor - say $F$ - of
 $M_f$ isomorphic to $M_n(\IC)$.  In the infinite algebra $N$ we find
 a pair of unital and commuting subfactors $B$ and $L$ of $N$ such
 that, $B$ is isomorphic to $B(\ell^2(\IN))$ and $N$ is isomorphic to
 the von Neumann algebra tensor product $B \bar{\otimes} L$.  Finally
 we will let $G$ denote the von Neumann subalgebra of $\pi(M) ^{\prime
   \prime}$ acting on K given as the sum $G = B \oplus F$. This is a
 von Neumann algebra of type I and hence injective. Further since $B$
 is infinite and $F$ is finite $Q$ is also a central projection in
 this algebra.  We can then average over the unitary translates
 $uhu^*$ of h with unitaries from $G$ and we can find a positive
 contraction $k$ in the commutant $G^{\prime}$ of $G$ such that

\begin{equation} \label{eq:5}
 k \in \overline{\text{conv}}^{uw}\{ uhu^* | u \,\mathrm{unitary in}
 \, G \} \cap G^{
 \prime}
 \end{equation}

First we remark that since $Q$ is a central projection in $G$
we must have that $k$ commutes with $Q$ and secondly we see from the
construction of $k$ and 
\ref{eq:3}, that this inequality must hold for $k$  too.

\begin{equation} \label{eq:6}
\forall u, v  \text{ unitaries in } \pi(M)^{ \prime \prime} ,
 \quad \quad  0 \leq uku^* \leq cvkv^*.
\end{equation}

The next observations with respect to $k$ have to be performed
according to the decomposition of the Hilbert space
 $K$ via $Q$ and $I-Q$, so we will split \ref{eq:6}.
Let us define $k_{\infty} =  k(I - Q)$ and $k_f = kQ$ then we get from
\ref{eq:6}
\begin{equation} \label{eq:7}
\forall u, v  \text{ unitaries in } N, 
 \quad \quad  0 \leq uk_{\infty}u^* \leq cvk_{\infty}v^*.
\end{equation}
\begin{equation} \label{eq:8}
\forall u, v  \text{ unitaries in } M_f ,
 \quad \quad  0 \leq uk_fu^* \leq cvk_fv^*.
\end{equation}

With respect to \ref{eq:7} it is quite easy to see that this one extends
to unitaries in any matrix algebras over $N$ with the same constant $c$
because  $k_{\infty}$ is in the commutant of $B$ and $B$ is an infinite
 tensor
factor of $N$. Hence we get immediately the following inequality with
respect to $\xi_{\infty}$.

\begin{equation} \label{eq:9}
  0 \leq ((\pi_{\infty})_{n}(V) \widetilde{k_{\infty}}(\pi_{\infty})_{n}(V^*)
  \xi_{\infty}, \xi_{\infty}) \leq
 (c\widetilde{k_{\infty}}\xi_{\infty},\xi_{\infty}). 
\end{equation}

We will now show that we can obtain a similar inequality with respect to
$\xi_f$. Here we will use an argument which is based on the one Pisier
uses in the proof of \cite[Lemma 5]{Pi7} where it is proved  
that an operator in $M_n(M)$
which is
supported on a projection of trace $1/n$  ( with respect to the normalized
trace ) can be factored in the way described in Definition 2.1.

We then define partial isometries $w_l, \,  1 \leq l \leq n$ in $M_n(M_f)$.
Let $(e_{ij})$ denote a set of matrix units for the $M_n(\IC)$ part
of the product $M_n(M) = M_n(\IC) \otimes M_f$, then we can write the the
partial isometries $w_l$ as sums of tensors as below and we get

\begin{displaymath}
w_l = \sum_{i=1}^n e_{i1}\otimes f_{li}; \quad w_lw_l^* = \widetilde{f_{ll}};
\quad w_l^*w_l = e_{11}\otimes I_{M_f}.
\end{displaymath}

From the inequalities \ref{eq:4} we get since $ p_j = f_{jj}$ and $ f_{jj}$
commutes with $k_f$ and $ k_f $ is a contraction  that

\begin{align} \label{eq:10} 
   & ((\pi_f)_n(V) \widetilde{k_f}(\pi_f)_n(V^*)\xi_f, \xi_f) & \\ 
 \notag  = &\sum_{i,j,l = 1}^n (\widetilde{f_{ii}}(\pi_f)_n(V) \widetilde{f_{jj}}
  \widetilde{k_f}(\pi_f)_n(V^*)\widetilde{f_{ll}}\xi_f, \xi_f) & \\ 
 \notag   \leq &
  \sum_{l = 1}^n (\widetilde{f_{ll}}(\pi_f)_n(V) \widetilde{f_{ll}}
      \widetilde{k_f}(\pi_f)_n(V^*)\widetilde{f_{ll}}\xi_f, \xi_f) +
      \varepsilon.&
\end{align}
Since $w_l$ commutes with $\widetilde{k_f}$ we can use the equality
 $\widetilde{f_{ll}} = w_lw_l^*$ to write 
$\widetilde{f_{ll}}\widetilde{k_f} = w_l\widetilde{k_f}w_l^*$ and then
 transform \ref{eq:10} into a set of
 $n$ inequalities regarding operators on $K_f $ by identifying $K_f$
 with the subspace of $\IC^n \otimes K_f $ corresponding to the projection
 $e_{11}\otimes I_{K_f}$. We will now look at each of the terms in the
 sum above. Hence we identify  define vectors

 $\eta_l = w_l^*\xi_f$ in $K_f$ and operators $x_l$ in $M_f$ by
 \begin{displaymath}
 x_l = w_l^*(\pi_f)_n(V)w_l = \sum_{i,j = 1}^n f_{il}\pi_f(v_{ij})f_{lj}
 \end{displaymath}
 By construction, each operator $x_l$ is a contraction and hence since
 $M_f$ is a finite von Neumann algebra there exist unitary operators - say
 - $ y_l , z_l $ in $M_f$ such that $x_l =\frac{1}{2}(y_l + z_l )$.

 The inequalities \ref{eq:10} and \ref{eq:8} then  yields

 \begin{align} \label{eq:11}
  & ((\pi_f)_n(V) \widetilde{k_f}(\pi_f)_n(V^*)\xi_f, \xi_f) \\
 \notag \leq & \sum_{l=1}^n (x_lk_fx_l^*\eta_l,\eta_l) + \varepsilon \\
 \notag \leq &  \frac{1}{4}\sum_{l=1}^n ((y_l+z_l)k_f(y_l+z_l)^*\eta_l,\eta_l)
+   \frac{1}{4}\sum_{l=1}^n ((y_l-z_l)k_f(y_l-z_l)^*\eta_l,\eta_l)
+ \varepsilon  \\
\notag = & \frac{1}{2}\sum_{l=1}^n (y_lk_fy_l^*\eta_l,\eta_l) +
 \frac{1}{2}\sum_{l=1}^n (z_lk_fz_l^*\eta_l,\eta_l) + \varepsilon \\
\notag \leq &  c\sum_{l=1}^n (k_f\eta_l,\eta_l) + \varepsilon \\
\notag = & c(\widetilde{k_f}\xi_f,\xi_f) +\varepsilon.
 \end{align}
The last equality is due to the fact that $k_f$ commutes with
all the $f_{ij}$.
 Going back to the inequality \ref{eq:3} we get $ h \leq ck$ and 
$k \leq ch$  so 
 \ref{eq:9}
 and \ref{eq:11} show
 that
 \begin{equation} \label{eq:12}
 \pi_n(V)\widetilde{h}\pi_n(V)^* \leq c^3\tilde{h}
 \end{equation}
 and consequently $\phi$ is completely bounded and satisfies $\| \phi \|_{cb}
 \leq c^{\frac{3}{2}} = \| \phi \|^3. $ Hence the similarity degree is 3 and
 the length is too.
\end{proof}

 \section{Von Neumann factors with property $\Gamma$ as complemented
   subspaces.} \label{sec:compl}

As mentioned in the introduction we here consider a continuous finite
factor $M$ with property $\Gamma$ which is a complemented subspace of a
C*-algebra $A$ and we will prove that there exists a completely
positive projection of norm one from $A$ onto $M$. The way we prove it, is by
showing that there exists a completely bounded projection from $A$
onto $M$ and then refer to the articles \cite{CS1,CS2,Pi2,Pi3} to get the
result. The completely bounded projection is obtained as a
point ultra weak limit of a bounded net of continuous 
 linear mappings of  $A$ into  $M$.

\begin{Theorem} \label{Thm:proj}
Let $A$ be a unital C*-algebra on a Hilbert space $H$ and let $M$ be
a continuous finite factor on $H$ which is a subalgebra of $A$ such
that there exists a bounded projection  of $A$ onto $M$. If $M$ has
property $\Gamma$ then
there exists a completely positive projection of norm one from $A$  
onto $M$.
\end{Theorem}
\begin{proof}
Let $\pi$ denote a bounded projection from  $A$ onto $M$.
We will first prove that we may assume that $M$ is generated - as a von
Neumann algebra - by a
countable set of operators.
Suppose  that the theorem  has been proven for such von Neumann
algebras. 
Then for a general finite
factor as $M$ and 
any von Neumann subalgebra
say $N$ of $M$ there exists a conditional expectation of norm one  $\rho_N$ of $M$
onto $N$ coming from the trace by 

$ \text{tr}(\rho_N(m)n) = \text{tr}(mn) $ for 
$ m \in M$ and $n \in N$,  so in particular for each subfactor
$N$ of $ M$ with a countable set of generators we have the  bounded
projection
$\rho_N \pi$ of $A$ onto $N$. Moreover it follows nearly 
immediately from the
construction of $\rho_N$ that the property $\Gamma$ of $M$ is
inherited by $N$. Then 
 by assumption there  exists a
projection of norm one from $A$ onto $N$ and since the unit ball in
$M$ is ultra weakly compact we can perform a limit over norm one
projections from $A$ onto larger and larger subalgebras $N $ of $M$
and by compactness obtain a net of
projections of norm one from $A$ onto subalgebras of $M$ which
converges pointwise ultra weakly  to a
projection of norm one from $A$ onto $M$.

Let us now assume that M is generated by a countable set of operators.
Then there exists a countable set $(m_i)_{(i \in \IN)}$  which is ultra
weakly dense in the unit ball of $M$. The result of Dixmier \cite{Di} then
implies 
\begin{align} \label{eq:31}
&\forall n \in \IN \, \exists \, \mathrm{projections} \,\, p^n_k, \quad 1 \leq
k \leq n \quad \mathrm{such  \, that} \\
\notag & \text{tr}(p_k^n) = \frac{1}{n},\quad p^n_kp^n_j = 0 \quad \mathrm{if}\, k \ne j, \quad
\sum_{k=1}^n p^n_k = I \quad \mathrm{ and } \\ 
\notag & \forall i , k  \in \{1, \dots , n\} \quad  \|
[p^n_k, m_i] \|_2 \leq \frac{1}{n^2}.
\end{align}
The choice of the projections $p^n_k$ and the density of the sequence
$(m_i)_{(i \in \IN)}$ imply  that we get the following result

\begin{equation} \label{eq:32}
\forall m \in M \quad  \sum_{k=1}^n p^n_kmp^n_k \quad \text{ converges
  ultrastrongly  to } m \text{ for }  n \to  \infty.
\end{equation} 

In order to be able to use \ref{eq:32} we have to make a modification
of $\pi$ corresponding to each of the sets $p^n_k, 1 \leq k \leq
n$. Hence we choose for each $n \in \IN$ a set of matrix units - say
$f^n_{ij}$ - for a subfactor - say $F^n$ - of  $M$ such that $F^n$ is
isomorphic to $M_n(\IC)$ and 
$\forall k \in \{ 1, \dots , n \}: f^n_{kk} = p^n_k.$ Let $U^n$
denote the group of unitaries in $F^n$, then this is a compact group
and it consequently has a Haar probability measure say $\mu$. we can
now define  a projection $\pi^n$ of $A$ onto $M$ which is modular with
respect to elements from $F^n$ by

\begin{equation} \label{eq:33}
\forall a \in  A \quad \pi^n(a) = \int_{ u \in U^n}\int_{ v \in U^n}
u\pi(u^*av)v^* \, d\mu(u)d\mu(v).
\end{equation} 

We can not prove that these projections are completely bounded but we
can construct a sequence  $\rho^n$ of mappings from $A$ into $M$ which has
the property that $\|\rho^n\|_n \leq \| \pi \| $ and $\forall m \in M
 \, \rho^n(m) \to m $ ultra strongly for $ n \to  \infty, n \in \IN$. We define
 $\rho^n$ by 

\begin{equation} \label{eq:34}
\forall a \in A \quad \rho^n(a) = \sum_{k=1}^n f^n_{kk}\pi^n(a)f^n_{kk}.
\end{equation}

and prove that $\|\rho^n\|_n \leq \| \pi \| $. In order to do so we
fix a set of matrix units $e_{jk}$ for $M_n(\IC)$ and 
define partial isomteries $w^n_k \in M_n(\IC) \otimes F^n $ by 
$w^n_k = \sum_{i=1}^n  e_{i1} \otimes f^n_{ki}$. Then $w^n_k{w^n_k}^* =
I \otimes f^n_{kk}$ and ${w^n_k}^*w^n_k = e_{11} \otimes I_M.$ In the
following computations we will identify the algebra
$(e_{11} \otimes I )[ M_n(\IC) \otimes M ] (e_{11} \otimes I )$ with  $M$
in the natural way and hence we get for  $X = (x_{ij})  \in M_n(A)$ 
\begin{equation} \label{eq:35}
\rho^n_n(X) = \sum_{k=1}^n(I \otimes f^n_{kk})\pi^n_n(X)(I \otimes
f^n_{kk}) 
\end{equation}
Hence, since we are working with the operator norm;
\begin{equation} \label{eq:36}
\|\rho^n_n(X)\| = \, \mathrm{max}\{ \|(I \otimes f^n_{kk})\pi^n_n(X)(I \otimes
f^n_{kk})\| \, | 1 \leq k \leq n \}.
\end{equation}

Now for each k we use the properties of $w^n_k$ and the $F^n$
modularity of $\pi^n$ to see   that 

\begin{equation} \label{eq:37}
\| (I \otimes f^n_{kk})\pi^n_n(X)(I \otimes
f^n_{kk}) \| = \| {w^n_k}^*\pi^n_n(X)w^n_k \| = \| \pi( \sum_{i,j = 1
  }^n f_{ik}x_{ij}f_{kj}) \|. 
\end{equation}

The last sum inside $\pi$ is obtained in $M$ via the identification
mentioned above and the norm of the sum is dominated by $\|X\|$ since
the sum is nothing but ${w^n_k}^*Xw^n_k$. A combination of this and
the results \ref{eq:36} and \ref{eq:37} give that

\begin{equation} \label{eq:38}
\| \rho^n_n \| \leq \| \pi \|
\end{equation}

The sequence of uniformly bounded mappings $\rho^n$ of $A$ into $M$
has a subnet which converges pointwise ultra weakly to a linear mapping
say $\rho$
 of $A$ into $M$. By \ref{eq:38}  we get that $\|\rho\|_n \leq
 \|\pi\|$, so $\rho$ is completely bounded. Further we get from
 \ref{eq:32} and the fact that the $\pi^n$ all are projections onto
 $M$ that $\forall m \in M \, \rho(m) = m $. We have then proved that
 $\rho$ is a completely bounded projection from $A$ onto $M$. By \cite{CS1,Pi2}
 it then follows that there exists a completely positive
 projection of norm one from $A$ onto $M$ and the proof is completed.
\end{proof}

 \section{Continuous Hochschild cohomology of von Neumann factors with property $\Gamma$.} \label{sec:coho}

The preprint \cite{CS3} which was never published contains the result
that for a von Neumann factor $M$ with property $\Gamma$ the second
continuous Hocschild 
cohomology group of  $M$ with coefficients in $M$,  $H^2_c(M,M)$
vanishes. The result was later published in the book \cite{SS1} by
Sinclair and Smith. We think that the methods used above should be
applicable as an ingredient in a proof of a general vanishing theorem
for the continuous Hocschild cohomology groups of factors with
property $\Gamma$. We are not able to get that far but we can get a
new and quite easy proof of the fact that  $H^2_c(M,M) = 0$ by using
the methods above and the nice result from \cite[ Proof of Theorem 5.1 ]{SS2} which in a very
short form says; that in order to
show that a continuous $n$-cocycle is a coboundary it is sufficient $(n-1)$
times to
prove that it is cohomologous to one which is completely bounded in
one variable only.

\medskip
  
\begin{Theorem} \label{Thm:h2}
Let M be a continuous finite von Neumann factor with property $\Gamma$ then 
the continuous Hochschild cohomolgy group $H^2_c(M,M)$ vanishes.
\end{Theorem}

\begin{proof}
  As in the case above where the algebra is a complemented subspace of
  some larger C*-algebra we will like to show first that it is
  sufficient to prove the result for a von Neumann factor which is
  countably generated.  

\smallskip So suppose that the result has been
  established in this case and let $\Phi$ denote a continuous
  2-cocycle on $M$. Then for any subalgebra $N$ of $M$ we have - as
  demonstrated in the proof of \ref{Thm:proj} - a completely positive
  projection $\pi_N$ of norm one from $M$ onto $N$.  This projection
  is also an $N$-bi-module mapping, so a simple algebraic manipulation
  shows that the composed bilinear map given by 
 $\pi_N\Phi : N \times  N \to N$ 
 is a continuous 2 cocycle on $N$. In the proof below we
  will show that when $N$ is countably generated then this cocycle is
  the coboundary of a continuous linear mapping $\phi_N : N \to N$
  which satisfies $\| \phi_N \| \leq 651\| \Phi \| $. If one defines
  $\psi_N = \phi_N\pi_N $ one gets a bounded net of mappings of $M$
  into $M$ indexed by the set of countably generated subalgebras of
  $M$. We may then find a subnet which converges pointwise ultra weakly
  to a continuous 1-cochain $\psi$ on $M$ and in turn get that $\Phi$ is
  the coboundary of $\psi$ 
  
  \medskip 
  
  Let us now suppose that $M$ is countably generated and that M acts
  standardly on $H$ with $\xi$ a cyclic and separating unit trace
  vector for M. The involution induced by $\xi$ is denoted $J$. 

The start of the proof follows the proof of  \cite[Theorem
  6.4.2]{SS1} where the stage is set.   
The pre dual
  of $M$ is now separable and by Popa's result \cite[Corollary 4.1]{Po} there
  exists an injective subfactor $R$ in $M$ such that 
$R^{\prime}\cap M   = \IC I_M$ 

By \cite[Theorem 3.1.1]{SS1} we may assume that $\Phi$ is
multimodular with respect to $R$, separately ultra weakly continuous
and vanishes whenever any of the arguments is in $R$. Further by
\cite[Corollary 5.2.4]{SS1} there is a an ultra weakly continuous
$R$-bimodular mapping $\phi :M \to JR^{\prime}J$ such that the coboundary
$\partial \phi$ equals $\Phi$. Further by \cite[Lemma 3.4.2]{SS1} we
can choose $\phi$ such that 

\begin{equation} \label{eq:50}
 \|\phi\| \leq 65 \|\Phi\| \quad \mathrm{ and \,\,trivially }\quad
 \|\Phi\| \leq 3\|\phi\|
\end{equation}

Since we are only dealing with 2-cocycles some direct estimates can be
made to show that much less than a factor of 65 will do as well.

By \cite[Lemma 5.4.7, (2)]{SS1} we get
a norm estimate for the action of 
$\phi_k : M_k(\IC) \otimes M \to M_k(\IC) \otimes JR^{\prime}J$, but in
order to understand the inequality below we must say that we use the
term $\| \,\|_2$ on any finite factor to mean the 2-norm with respect
to the trace state - or the normalized trace on the algebra. Then we
can quote \cite{SS1} as:
 
\begin{equation} \label{eq:51}
\forall  k \in \IN \quad \forall x \in M_k(\IC) \otimes M : \quad
\|\phi_k(x)\|^2 \leq 4\| \phi \|^2( \|x\|^2 + k\|x\|^2_2).
\end{equation}  

We are now in the position to use the $\Gamma$-property in a similar
way as above. Let $\{ m_i | i \in \IN \}$ be a sequence in the unit
ball of $M$ which is dense with respect to the $\| \, \|_2 $ topology.
For each $n  \in \IN$ we choose using \cite{Di} 
a set of pairwise orthogonal and
equivalent projections $\{ p^n_1, \dots , p^n_n \}$ in $M$ 
with sum $I$  such that

\begin{equation} \label{eq:52}
\forall i, j, k  \in \{ 1, 2, \dots , n \} \quad \| [m_j, p_i^n] \|_2 \leq
 n^{-3} \quad \mathrm{and} \quad \quad \| [\Phi(m_j, m_k), p^n_i] \|_2 \leq
 n^{-3}.
\end{equation}

For each $n \in \IN$ we will modify $\Phi$ by a coboundary  say
$\Psi^n$ which is related to the set of projections 
$\{ p^n_1, \dots , p^n_n \}$ in such  a way that we can prove that $\Phi $ is
cohomologous to a 2-cocycle which is completely bounded in the first
variable. By \cite{SS2} this is  sufficient in order to see  that $\Phi$
is a coboundary too.

\begin{equation} \label{eq:53}
\forall n \in \IN \, \forall m \in M \quad \mathrm{let} \quad \psi^n(m) =
\sum_{i=1}^n p^n_i\Phi(p^n_i, m)p^n_i.
\end{equation}

This $\psi^n$ is clearly an ultra weakly  continuous linear map of $M$
into $M$ such that $\| \psi^n \| \leq \| \Phi \|$, consequently for the
coboundary - say $ \Psi^n = \partial \psi^n$ -  we have  
$\| \Psi^n \| \leq 3 \| \Phi \|$.

 In order to clarify
the following computations we introduce some 
 bilinear operators
from $ M \times M $ to $M$ by 

\begin{align} \label{eq:54}
\Omega^n & = \Phi - \Psi^n \\
\notag \Delta^n_1(x,y) & = \sum_{i=1}^n p^n_i\Phi(p^n_i, x)[y, p^n_i] \\
\notag \Delta^n_2(x,y) &=  \Phi(x,y) - \sum_{i=1}^n p^n_i\Phi(x,y)
p^n_i \\
\notag \Delta^n & =  \Delta^n_1 + \Delta^n_2 \\
\notag \Theta^n(x,y) & = \sum_{i=1}^n (p^n_i\Phi(p^n_ix, y)p^n_i - 
xp^n_i\Phi(p^n_i,y)p^n_i.
\end{align}

We can then examine $\Omega^n$. The fact that $\Phi $ is a
2-cocycle is used from the second to the third line just below.

\begin{align} \label{eq:55}
\Omega^n(x, y) & = \Phi(x, y) + \sum_{i=1}^n(-xp_i^n\Phi(p_i^n,
y)p_i^n + p_i^n\Phi(p_i^n, xy)p_i^n - p_i^n\Phi(p_i^n, x)p_i^n y) \\
\notag & = \Phi(x, y) + \Delta^n_1(x, y)+ \sum_{i=1}^n(-xp_i^n\Phi(p_i^n,
y)p_i^n + p_i^n\Phi(p_i^n, xy)p_i^n - p_i^n\Phi(p_i^n, x)yp_i^n )\quad
  \\
\notag & = \Phi(x, y)+ \Delta^n_1(x, y)+ \sum_{i=1}^n(-xp_i^n\Phi(p_i^n,
y)p_i^n - p_i^n\Phi(x, y)p_i^n + p_i^n\Phi(p_i^nx, y)p_i^n) \\
\notag & = \Delta^n(x,y) + \Theta^n(x, y).
\end{align} 

We know that $\Phi  = \partial \phi$ so we can see that $\Theta^n $ is
expressed by 

\begin{align} \label{eq:56}
\Theta^n(x,y)& = \sum_{i=1}^n (p^n_ix\phi(y)p^n_i
-p^n_i\phi(p^n_ixy)p^n_i + p^n_i\phi(p^n_ix)yp^n_i 
- xp^n_i\Phi(p^n_i,y)p^n_i) \\
\notag & =  \sum_{i=1}^n (p^n_ix\phi(y)p^n_i -
xp^n_i\Phi(p^n_i,y)p^n_i) + \sum_{i=1}^np^n_i(-\phi(p^n_ixy) +
\phi(p^n_ix)y)p^n_i.
 \end{align}

This decomposition shows that for a fixed y in $M$ and a fixed natural
number $k$ we get using  \ref{eq:51} and \ref{eq:56} that 
for 
$X \in M_k(\IC)\otimes M, \quad  \widetilde{y} = I_{M_k(\IC)} \otimes y $ 
and any $  n \geq k \quad   n \in \IN$

\begin{align} \label{eq:57} 
\|\Theta^n_k(X,\widetilde{y})\| & \leq 4\|\phi\|\|X\|\|y\| +
\mathrm{max}\{\, \|\widetilde{p^n_i} \phi_k(\widetilde{p^n_i}X\widetilde{y})\widetilde{p^n_i}\| +
\| \widetilde{p^n_i}\phi_k(\widetilde{p^n_i}X)\widetilde{y}\widetilde{p^n_i}\|
| 1 \leq i \leq n\} \\
\notag & \leq  4\|\phi\|\|X\|\|y\| + 2(2\sqrt{2})\|\phi\|\|X\|\|y\| \\
\notag & \leq 10\|\phi\|\|X\|\|y\|.
\end{align}

It is clear from the construction of $\Delta^n$ in \ref{eq:54} that this
sequence of bilinear operators on $M$ is uniformly bounded and
converges pointwise ultra strongly towards $0$. A combination of this
with the estimates from \ref{eq:57} shows that if we take a subnet of
the sequence $(\psi^n)_{n \in \IN}$ which converges pointwise
ultra weakly to a continuous 1-cochain - say $\psi $ on $M$ then the
2-cocycle - say $\Omega$ on $M$  given by
$\Omega = \Phi - \partial \psi $ is completely bounded in the left variable,
and by the methods from \cite[Proof of Theorem 5.1]{SS2} we can show
that a 2-cocycle which is completely bounded in the left variable is
the coboundary of a continuous 1-cochain - say  $\omega $ - 
on $M$  which is in
 the pointwise ultra weakly closed convex hull of the set of continuous
 1-cochains of  the form below;
 
\begin{equation} \label{eq:58}
 m \to \sum_{j \in J}x_j^*\Omega(x_j, m) \quad \mathrm{where} \quad x_j \in M
 \quad \mathrm{and} \quad \sum_{j \in J}x_j^*x_j = I_M.
\end{equation}

We do have $\Phi = \partial(\psi + \omega) $ and by a  combination of
\ref{eq:50}, \ref{eq:53}, \ref{eq:57} and  \ref{eq:58} we get that  

\begin{align} \label{eq:59}
\|\psi + \omega \| & \leq \|\Phi\| + \|\omega\| \\
\notag & \leq \|\Phi\|+ 10\|\phi\| \\
\notag & \leq \|\Phi\| + 650\|\Phi\| \\
\notag & = 651\|\Phi\
\end{align}  

We have now proved that if $ M $ is countably generated any continuous
2-cocycle $\Phi$ is inner and we have moreover obtained a
universal bound on the cochains implementing $\Phi$ so we may conclude
that the theorem is proved for a general continuous von Neumann factor
with property $\Gamma$.

\end{proof}

\end{document}